\documentclass[11pt]{amsart}
\usepackage[utf8]{inputenc}
\usepackage{enumerate,amssymb,xcolor,comment}
\usepackage[a4paper, footskip=1cm, headheight = 16pt, top=3.5cm, bottom=2.5cm,  right=2.5cm,  left=2.5cm ]{geometry}
\usepackage[colorlinks,linkcolor=blue,citecolor=red]{hyperref}
\usepackage[center]{caption}

\newtheorem{theorem}{Theorem}[section]
\newtheorem{lemma}[theorem]{Lemma}
\newtheorem{conjecture}[theorem]{Conjecture}

\usepackage[foot]{amsaddr}

\usepackage[T1]{fontenc}

\usepackage[leqno]{amsmath}

\makeatletter
\newcommand{\leqnomode}{\tagsleft@true}
\newcommand{\reqnomode}{\tagsleft@false}
\makeatother

\usepackage{latexsym}
\usepackage{amsfonts}

\usepackage{tikz}
\usepackage{float}

\def\dd{\hbox{-}}

\usepackage{marvosym}               

\newcounter{tbox}

\newcommand{\sta}[1]{\medskip\medskip\refstepcounter{tbox}\noindent{\parbox{\textwidth}{(\thetbox) \emph{#1}}}\vspace*{0.3cm}}

\newcommand{\mylongtitle}[1]{%
  \ifodd\value{page}%
    \protect\parbox{0.97\linewidth}{#1}\hfill%
  \else%
    \hfill\protect\parbox{0.97\linewidth}{#1}%
  \fi%
}

\usepackage{latexsym}
\usepackage{amsfonts}
\usepackage[leqno]{amsmath}
\usepackage{tikz}
\usepackage{float}
\usepackage{lmodern}
\usepackage[T1]{fontenc}

\def\dd{\hbox{-}}

\usepackage{marvosym}

\newcommand{\otherlabel}[2]{\protected@edef\@currentlabel{#2}\label{#1}}
\makeatother

\mathchardef\mh="2D

\title[Hitting all maximum stable sets in $P_5$-free graphs]{Hitting all maximum stable sets in $P_5$-free graphs\footnote{This is an accepted manuscript. The published version appeared in Journal of Combinatorial Theory, Series B, 
Volume 165, March 2024, Pages 142-163, \url{https://doi.org/10.1016/j.jctb.2023.11.005}.}}

\author{Sepehr Hajebi $^{\ast}$}
\author{Yanjia Li $^{\ast\ast}$}
\author{Sophie Spirkl$^{\ast \dagger}$}

\address{$^{\ast}$Department of Combinatorics and Optimization, University of Waterloo, Waterloo, Ontario, Canada}
\address{$^{\ast\ast}$School of Mathematics,
Georgia Institute of Technology,
Atlanta, GA 30332, USA}
\address{$^{\dagger}$ We acknowledge the support of the Natural Sciences and Engineering Research Council of Canada (NSERC), [funding reference number RGPIN-2020-03912].
Cette recherche a \'et\'e financ\'ee par le Conseil de recherches en sciences naturelles et en g\'enie du Canada (CRSNG), [num\'ero de r\'ef\'erence RGPIN-2020-03912]. This project was funded in part by the Government of Ontario.}
\date {\today}

\begin{document}

\maketitle
\begin{abstract}
We prove that every $P_5$-free graph of bounded clique number contains a small hitting set of all its maximum stable sets (where $P_t$ denotes the $t$-vertex path, and for graphs $G,H$, we say $G$ is \textit{$H$-free} if no induced subgraph of $G$ is isomorphic to $H$).

More generally, let us say a class $\mathcal{C}$ of graphs is \textit{$\eta$-bounded} if there exists a function $h:\mathbb{N}\rightarrow \mathbb{N}$ such that $\eta(G)\leq h(\omega(G))$ for every graph $G\in \mathcal{C}$, where $\eta(G)$ denotes smallest cardinality of a hitting set of all maximum stable sets in $G$, and $\omega(G)$ is the clique number of $G$. Also, $\mathcal{C}$ is said to be \textit{polynomially $\eta$-bounded} if in addition $h$ can be chosen to be a polynomial.

We introduce $\eta$-boundedness inspired by a question of Alon (asking how large $\eta(G)$ can be for a $3$-colourable graph $G$), and motivated by a number of meaningful similarities to $\chi$-boundedness, namely,
\begin{itemize}
\item given a graph $G$, we have  $\eta(H)\leq \omega(H)$ for every induced subgraph $H$ of $G$ if and only if $G$ is perfect;
    \item there are graphs $G$ with both $\eta(G)$ and the girth of $G$ arbitrarily large; and
    \item if $\mathcal{C}$ is a hereditary class of graphs which is polynomially $\eta$-bounded, then $\mathcal{C}$ satisfies the Erd\H{o}s-Hajnal conjecture.
\end{itemize} 
The second bullet above in particular suggests an analogue of the Gy\'{a}rf\'{a}s-Sumner conjecture, that the class of all $H$-free graphs is $\eta$-bounded if (and only if) $H$ is a forest. Like $\chi$-boundedness, the case where $H$ is a star is easy to verify, and we prove two non-trivial extensions of this: $H$-free graphs are $\eta$-bounded if (1) $H$ has a vertex incident with all edges of $H$, or (2) $H$ can be obtained from a star by subdividing at most one edge, exactly once.

Unlike $\chi$-boundedness, the case where $H$ is a path is surprisingly hard. Our main result mentioned at the beginning shows that $P_5$-free graphs are $\eta$-bounded. The proof is rather involved compared to the classical ``Gy\'{a}rf\'{a}s path'' argument which establishes, for all $t$, the $\chi$-boundedness of $P_t$-free graphs. It remains open whether $P_t$-free graphs are $\eta$-bounded for $t\geq 6$.

It also remains open whether $P_5$-free graphs are polynomially $\eta$-bounded, which, if true, would imply the Erd\H{o}s-Hajnal conjecture for $P_5$-free graphs. But we prove that $H$-free graphs are polynomially $\eta$-bounded if $H$ is a proper induced subgraph of $P_5$. We further generalize the case where $H$ is a $1$-regular graph on four vertices, showing that $H$-free graphs are polynomially $\eta$-bounded if $H$ is a forest with no vertex of degree more than one and at most four vertices of degree one.

\end{abstract}


\section{Introduction}\label{intro}
\subsection{Background and motivation} Throughout this paper, all graphs have finite vertex sets and no loops or parallel edges. We denote by $\mathbb{N}$ the set of all positive integers. Let $G=(V(G),E(G))$ be a graph. A \textit{stable set} in $G$ is a set of pairwise non-adjacent vertices. The cardinality of the largest stable set in $G$ is denoted by $\alpha(G)$, and a \textit{maximum stable set} in $G$ is a stable set in $G$ of cardinality $\alpha(G)$. By a \textit{hitting set} for $G$, we mean a set $X\subseteq V(G)$ which intersects every maximum stable set in $G$, or equiavalenty, satisfies $\alpha(G\setminus X)<\alpha(G)$. We denote by $\eta(G)$ the smallest cardinality of a hitting set for $G$. Bollob\'{a}s, Erd\H{o}s and Tuza raised the following conjecture, that in any graph $G$ with a stable set which is larger than $V(G)$ by at most a constant factor, there exists a hitting set which is much smaller than the entire vertex set.

\begin{conjecture}[Bollob\'{a}s, Erd\H{o}s and Tuza \cite{Erdlegacy,Erdproblems}]\label{Erdosconj}
For every real $\delta>0$, every graph $G$ with $\alpha(G)\geq \delta|V(G)|$ satisfies $\eta(G)\leq o(|V(G)|)$.
\end{conjecture}

This is still open, and might be true as far as we know. In fact, a result of Hajnal \cite{Hajnal} (see also \cite{Rabern}) shows that Conjecture~\ref{Erdosconj} holds for all $\delta> 1/2$. Also, the following recent result of Alon \cite{Alon} complements Conjecture~\ref{Erdosconj}, providing a construction of infinitely many graphs $G$ with $\alpha(G)$ linear in $|V(G)|$ and $\eta(G)$ arbitrarily large.

\begin{theorem}[Alon \cite{Alon}]\label{Alonlowerbound}
For every integer $n\geq 1$, there exists a graph $G_n$ with $|V(G_n)|>n$, $|V(G_n)|/2 >\alpha(G_n)> |V(G_n)|/4$ and $\eta(G_n)>\sqrt{|V(G_n)|}/2$.
\end{theorem}

One approach to understand Conjecture~\ref{Erdosconj} is to examine it for special classes of graphs with $\alpha$ linear in the number of vertices, where the most natural candidates are graph classes of bounded chromatic number (the chromatic number of a graph $G$ is denoted by $\chi(G)$; recall that every graph $G$ satisfies $\alpha(G)\geq |V(G)|/\chi(G)$). As a starter, Hajnal's result \cite{Hajnal} mentioned above implies that every graph $G$ with $\chi(G)\leq 2$ satisfies $\eta(G)\leq 2$. But the train already stops here; as highlighted by Alon \cite{Alonober}, it is not known whether $3$-colourable graphs satisfy Conjecture~\ref{Erdosconj}. Although bounded $\eta$ is not expected anymore, as Alon~\cite{Alonprivate} shared with us, in private communication, a proof of the following result:
\begin{theorem}[Alon \cite{Alonprivate}]
    For every integer $h\geq 1$, there exists an $8$-colourable graph $G_h$ with $\eta(G_h)> h$.
\end{theorem}

Note that the complete graph $K_t$ satisfies $\chi(K_t)=\eta(K_t)=t$. Inspired by $\chi$-boundedness as the study of graph classes in which complete graphs are the only (induced subgraph) obstructions to bounded chromatic number, we look into classes in which $\eta$ can be upper bounded by a function of the clique number.

\subsection{$\eta$-boundedness vs. $\chi$-boundedness}

\sloppy We begin with a brief digression about $\chi$-boundedness. For graphs $G$ and $H$, we say $H$ is an \textit{induced subgraph} of $G$ if $H$ can be obtained from $G$ by removing vertices, and that $G$ is \textit{$H$-free} if no induced subgraph of $G$ is isomorphic to $H$. A class $\mathcal{C}$ of graphs is said to be \textit{hereditary} if $\mathcal{C}$ is  closed under isomorphism and taking induced subgraphs. A \textit{clique} in a graph $G$ is a set of pairwise adjacent vertices. The cardinality of the largest clique in $G$ is denoted by $\omega(G)$, and a \textit{maximum clique} in $G$ is a clique in $G$ of cardinality $\omega(G)$.  A class $\mathcal{C}$ of graphs is said to be \textit{$\chi$-bounded} if there exists a function $f:\mathbb{N}\rightarrow \mathbb{N}$ such that every graph $G\in \mathcal{C}$ satisfies $\chi(G)\leq f(\omega(G))$. This is a heavily studied notion of interest in structural and extremal graph theory (see \cite{chibnddsurvey} for a survey), mainly due to the following three reasons.

\begin{itemize}
    \item First, $\chi$-boundedness has its roots in the theory of perfect graphs (that is, graphs $G$ with $\chi(H)\leq \omega(H)$ for every induced subgraph $H$ of $G$). In fact, the systematic study of $\chi$-boundedness was initiated in a 1987 paper by Gy\'{a}rf\'{a}s \cite{Gyarfaschi} entitled ``Problems from the world surrounding perfect graphs''.
    \item Second, it has an intriguing connection to a well-known conjecture of Erd\H{o}s and Hajnal. Let us say that a hereditary class $\mathcal{C}$ of graphs has the \textit{Erd\H{o}s-Hajnal property} if there exists $\varepsilon>0$ such that every graph $G\in \mathcal{C}$ satisfies either $\alpha(G)\geq |V(G)|^{\varepsilon}$ or $\omega(G)\geq |V(G)|^{\varepsilon}$. While the class of all graph is known not to have the Erd\H{o}s-Hajnal property, the Erd\H{o}s-Hajnal conjecture \cite{EH1,EH2} asserts that \textit{all} other hereditary classes have the Erd\H{o}s-Hajnal property, or equivalently, for every graph $H$, the class of all $H$-free graphs has the Erd\H{o}s-Hajnal property. This remains out of reach even when $H$ is the five-vertex path $P_5$. Back to the connection to $\chi$-boundedness, one may observe that a hereditary class of graphs $\mathcal{C}$ has the Erd\H{o}s-Hajnal property if $\mathcal{C}$ is \textit{polynomially-$\chi$-bounded}, that is, $\mathcal{C}$ is $\chi$-bounded where the function $f:\mathbb{N}\rightarrow \mathbb{N}$ from the definition can be chosen as $f(x)=x^d$ for some positive integer $d$. As expected, it is not known whether $P_5$-free graphs are polynomially-$\chi$-bounded.
    \item Third, despite the tremendous advancement in $\chi$-boundedness over the past 10 years, there remains an open problem at its heart, the Gy\'{a}rf\'{a}s-Sumner conjecture \cite{GGS,GSS}, asserting that the class of all $H$-free graphs is $\chi$-bounded if $H$ is a forest. This is easy to verify when $H$ is the star $K_{1,t}$ (on $t+1$ vertices), and was proved by Gy\'{a}rf\'{a}s \cite{Gyarfaschi} when $H$ is the path $P_t$ (on $t$ vertices) for all $t$. Also, the Gy\'{a}rf\'{a}s-Sumner conjecture would be sharp in the sense that forests are the only graphs that may possibly satisfy this conjecture; this is evidenced by a well-known result of Erd\H{o}s \cite{Erdprobability} that there are graphs with arbitrarily large girth and arbitrarily large chromatic number.
\end{itemize}
Let us now define ``$\eta$-boundedness'' formally. We say a class $\mathcal{C}$ of graphs is \textit{$\eta$-bounded} if there exists a function $f:\mathbb{N}\rightarrow \mathbb{N}$ such that every graph $G\in \mathcal{C}$ satisfies $\eta(G)\leq f(\omega(G))$, and \textit{polynomially $\eta$-bounded} if in addition we have $f(x)=x^d$ for some positive integer $d$. Our interest in $\eta$-boundedness is sparked by the fact that it shares with $\chi$-boundedness the three reasons mentioned above for being a compelling notion.

First, the class of perfect graphs can be viewed as the point of origin for $\eta$-boundedness, more or less in the same way it is so for $\chi$-boundedness. To see this, we need an easy lemma, which we leave to the reader to prove.

\begin{lemma}\label{etatimes}
    Let $d,h\in \mathbb{N}$ and let $G$ be a graph such that $\eta(G')\leq h$ for every induced subgraph $G'$ of $G$. Then there exists 
    $D\subseteq V(G)$ with $|D|\leq dh$ such that $\alpha(G\setminus D)\leq \alpha(G)-d$. In particular, we have $\alpha(G)\geq |V(G)|/h$.
\end{lemma}

\begin{theorem}\label{perfectthm}
    Let $G$ be a graph. Then $G$ is perfect if and only if we have $\eta(H)\leq \omega(H)$ for every induced subgraph $H$ of $G$.
\end{theorem}

\begin{proof}
    Lov\'{a}sz \cite{lovasz} proved that a graph $G$ is perfect if and only of if we have $\alpha(H)\omega(H)\geq |V(H)|$ for every induced subgraph $H$ of $G$. This in particular implies that a graph $G$ is perfect if and only if the complement of $G$ perfect.

    Now, assume that $G$ is a perfect graph and $H$ is a non-null induced subgraph of $G$; then $H$ is perfect, too. By  Lov\'{a}sz's result, the complement of $H$ is perfect. It follows that $V(H)$ can be partitioned into cliques $C_1,\ldots, C_s$ where $\alpha(H)=s\geq 1$. Thus, for all $i\in \{1,\ldots, s\}$, every maximum stable set of $H$ intersects $C_i$ in exactly one vertex. But then $C_1$ is hitting set for $H$, and so $\eta(H)\leq |C_1|\leq \omega(H)$.

    Conversely, assume that $G$ is a graph with  $\eta(H)\leq \omega(H)$ for every induced subgraph $H$ of $G$. Let $H$ be an arbitrary induced subgraph of $G$. It follows that $\eta(H')\leq \omega(H')\leq \omega(H)$ for every induced subgraph $H'$ of $H$. Thus, by Lemma~\ref{etatimes} applied to $G=H$ and $h=\omega(H)$, we have $\alpha(H)\omega(H)\geq |V(H)|$. But then by Lov\' asz's Theorem, $G$ is perfect. This completes the proof of Theorem~\ref{perfectthm}.
\end{proof}

Our second motivation for introducing $\eta$-boundedness is its connection to the Erd\H{o}s-Hajnal conjecture, which is identical to that of $\chi$-boundedness.

\begin{theorem}\label{etatoEH}
    Let $\mathcal{C}$ be a hereditary class of graph which is polynomially $\eta$-bounded. Then $\mathcal{C}$ satisfies the Erd\H{o}s-Hajnal conjecture.
\end{theorem}
\begin{proof}
   Since $\mathcal{C}$ is polynomially $\eta$-bounded, there exists $d\in \mathbb{N}$ such that every graph $G\in \mathcal{C}$ satisfies $\eta(G)\leq \omega(G)^d$. We claim that for $\varepsilon=(d+1)^{-1}$, every graph $G\in \mathcal{C}$ satisfies either $\alpha(G)\geq |V(G)|^{\varepsilon}$ or $\omega(G)\geq |V(G)|^{\varepsilon}$. To see this,  note that since $\mathcal{C}$ is hereditary, for every graph $G\in \mathcal{C}$ and every induced subgraph $G'$ of $G$, we have $\eta(G')\leq \omega(G')^d\leq \omega(G)^d$. Thus, by Lemma~\ref{etatimes} applied to $G$ and $h=\omega(G)^d$, we have $\alpha(G)\omega(G)^d\geq |V(G)|$, which in turn implies that either $\alpha(G)\geq |V(G)|^{(d+1)^{-1}}$ or $\omega(G)\geq |V(G)|^{(d+1)^{-1}}$, as desired.
\end{proof}

The third reason why we find $\eta$-boundedness interesting is a direct analogue of the Gy\'{a}rf\'{a}s-Sumner conjecture, attempting to characterize all graphs $H$ for which $H$-free graphs are $\eta$-bounded. The result below shows that this holds only if $H$ is a forest. This can be proved using a direct adaption of Erd\H{o}s' proof \cite{Erdprobability} for the existence of graphs with large girth and large chromatic number, which we omit.

\begin{theorem}\label{largeglargeeta}
    For every $h\in \mathbb{N}$, there exists a graph $G_{h}$ with girth more than $h$ and with $\eta(G_{h})> h$.
\end{theorem}

In view of Theorem~\ref{largeglargeeta}, we propose the following analogue of the  Gy\'{a}rf\'{a}s-Sumner conjecture for $\eta$-boundedness.

\begin{conjecture}\label{gseta}
    For every forest $H$, there exists a function $f:\mathbb{N}\rightarrow \mathbb{N}$ such that every $H$-free graph $G$ satisfies $\eta(G)\leq f(\omega(G))$.
\end{conjecture}

The main content of this paper is the proof of Conjecture~\ref{gseta} for an assortment of special forests. Like the original Gy\'{a}rf\'{a}s-Sumner conjecture, Conjecture~\ref{gseta} is easily seen to hold when $H$ is star, which we prove next. Let us first mention the following quantified version of Ramsey's theorem. A similar result with a slightly different bound and an almost identical proof has also appeared in \cite{polychi2}; we include a proof though for the sake of completeness.

\begin{theorem}\label{classicalramsey}
    For all $c,s\in \mathbb{N}$, every graph $G$ on at least $c^s$ vertices contains either a clique of cardinality $c$ or a stable set of cardinality $s$.
\end{theorem}
\begin{proof}
The proof is by induction on $s$ for fixed $c$. The cases $c=1$ and $s=1$ are easily seen to hold. So we may assume that $c,s>1$. Let $G$ be a graph on at least $c^s$ vertices with no clique of cardinality  $c$ and no stable set of cardinality $s$. Let $K$ be a maximum clique in  $G$; thus, we have $|K|\leq c-1$. For each $x\in K$, let $M_x$  be the set of all vertices in $G$ which are non-adjacent to $x$. It follows from the maximality of $K$ that $V(G)=\bigcup_{x\in K}(M_x\cup \{x\})$. Now, for every $x\in K$, $G[M_x]$ contains no clique of cardinality $c$ (as neither does $G$) and no stable set of cardinality $s-1$ (or otherwise $M_x\cup \{x\}$, and so $G$, contains a stable set of cardinality $t$). Consequently, by the induction hypothesis, we have $|M_x|< c^{s-1}$. But then $|V(G)|\leq (c-1)c^{s-1}<c^s$, a contradiction. This proves Theorem~\ref{classicalramsey}.
\end{proof}

Note also that for every vertex $v$ in every graph $G$, every maximum stable set of $G$ contains either $v$ or a neighbour of $v$. In other words:

\begin{theorem}\label{bndddegree0}
    Let $G$ be a graph and $v\in V(G)$ be a vertex of degree at most $d$ in $G$. Then we have $\eta(G)\leq d+1$.
\end{theorem}
Combined with Theorem~\ref{classicalramsey}, this immediately yields Conjecture~\ref{gseta} for stars:
\begin{theorem}\label{bndddegree}
    Let $s\in \mathbb{N}$ and $G$ be a $K_{1,s}$-free graph. Then every vertex of $G$ has degree less than $\omega(G)^s$. Consequently,  we have $\eta(G)\leq \omega(G)^s$.
\end{theorem}
We remark that Theorem~\ref{bndddegree0} also implies that every graph class in which the degeneracy is bounded by some function of the clique number is $\eta$-bounded. Proper minor-closed classes are well-known to have this property. Other examples include, for every positive integer $s$ and every forest $H$, the class of all $K_{s,s}$-free $H$-free graphs \cite{polychi1} (in fact, the main result of \cite{polychi1} yields polynomial $\eta$-boundedness in this case), and for every positive integer $s$ and every graph $H$,  the class of all $K_{s,s}$-free graphs with no induced subgraph isomorphic to a subdivision of $H$ \cite{KuhnOsthus}. The latter in particular includes the famous class of ``even-hole-free'' graphs. Another another result from \cite{EHBS} shows that every even-hole-free graph $G$ has a vertex of degree at most $2\omega(G)-2$, which along with Theorem~\ref{bndddegree0} implies that $\eta(G)\leq 2\omega(G)-1$.

But this concludes the list of forests for which Conjecture~\ref{gseta} can be proved easily: for all forests apart from  induced subgraphs of stars  (for which Conjecture~\ref{gseta} follows from Theorem~\ref{bndddegree}) and induced subgraphs of $P_4$ (for which Conjecture~\ref{gseta} follows from Theorem~\ref{perfectthm}, as $P_4$-free graphs are perfect), Conjecture~\ref{gseta} seems to be much harder to prove than the Gy\'{a}rf\'{a}s-Sumner conjecture. The case of paths on five or more vertices is particularly tempting as there is a short and elementary proof of the Gy\'{a}rf\'{a}s-Sumner conjecture for all paths \cite{Gyarfaschi}.

Our main result verifies Conjecture~\ref{gseta} when $H$ is the five-vertex path:

\begin{theorem}\label{mainP5}
    The class of $P_5$-free graphs is $\eta$-bounded.
\end{theorem}

The proof of Theorem~\ref{mainP5} is substantially longer and harder than its $\chi$-boundedness counterpart. It remains open whether $P_t$-free graphs are $\eta$-bounded for $t\geq 6$.

\begin{conjecture}\label{Ptconj}
    For every integer $t\geq 6$, the class of $P_t$-free graphs is $\eta$-bounded.
\end{conjecture}

Similarly, in the case of disconnected forests, the Gy\'{a}rf\'{a}s-Sumner conjecture appears to be more approachable than Conjecture~\ref{gseta}. It is straightforward to show that a forest $H$ satisfies the Gy\'{a}rf\'{a}s-Sumner conjecture if and only if all its components do so. We do not know whether the corresponding statement holds for Conjecture~\ref{gseta}; in fact, it remains open whether the disjoint union of two stars satisfies Conjecture~\ref{gseta}. But we can prove the following weakening:

\begin{theorem}\label{mainKrs1}
   Let $H$ be a graph in which some vertex is incident with all edges. Then $H$-free graphs are $\eta$-bounded.
\end{theorem}
Moreover, the Gy\'{a}rf\'{a}s-Sumner conjecture is known to hold for all trees of radius two \cite{KP}. In contrast, Conjecture~\ref{gseta} remains open even for trees of diameter three, but we can come close:
\begin{theorem}\label{mainFt1}
   Let $H$ be a tree which is obtained from a star by subdividing at most one edge, exactly once. Then $H$-free graphs are $\eta$-bounded.
\end{theorem}

Finally, concerning the connection to the Erd\H os-Hajnal conjecture, we would like to propose the following strengthening of Conjecture~\ref{gseta}. We remark that the corresponding extension of the Gy\'{a}rf\'{a}s-Sumner conjecture remains open, too (although it is known that there are $\chi$-bounded graph classes which are not polynomially $\chi$-bounded \cite{daviespoly}).

\begin{conjecture}\label{gsetapoly}
    For every forest $H$, there exists $d\in \mathbb{N}$ such that every $H$-free graph $G$ satisfies $\eta(G)\leq \omega(G)^d$.
\end{conjecture}

Note that by Theorems~\ref{perfectthm} and \ref{bndddegree}, Conjecture~\ref{gsetapoly} holds when $H$ is either a star or an induced subgraph of $P_4$. The case where $H$ is the five-vertex path remains open;  this is particularly interesting because by Theorem~\ref{etatoEH}, it would imply the Erd\H os-Hajnal conjecture for $P_5$-free graphs. Let us state it separately:

\begin{conjecture}\label{gsetapolyP5}
    There exists $d\in \mathbb{N}$ such that every $P_5$-free graph $G$ satisfies $\eta(G)\leq \omega(G)^d$.
\end{conjecture}

Indeed, the bound obtained in Theorem~\ref{mainP5} is super-exponential, and it would be interesting even to bring it down to a singly exponential bound.

On the bright side, we were able to prove the following:

\begin{theorem}\label{properP51}
    Let $H$ be isomorphic to a proper induced subgraph of $P_5$. Then $H$-free graphs are polynomially $\eta$-bounded.
\end{theorem}

For each $t\in \mathbb{N}$, let $M_t$ be the unique (up to isomorphism) $1$-regular graph on $2t$ vertices. Then $M_2$ is isomorphic to a proper induced subgraph of $P_5$, and so by Theorem~\ref{properP51}, $M_2$-free graphs are polynomially $\eta$-bounded. In this case, we can prove more:

\begin{theorem}\label{mainparallel1}
    Let $H$ be a graph with no vertex of degree more than one and at most four vertices of degree one. Then $H$-free graphs are polynomially $\eta$-bounded.
\end{theorem}

We would also like to conjecture a natural strengthening of Theorem~\ref{mainparallel1}, that $M_t$-free graphs are polynomially $\eta$-bounded for all $t$, though it is not even known (stikingly enough) whether $M_3$-free graphs are $\eta$-bounded.

\begin{conjecture}\label{Mtconj}
    For every integer $t\geq 3$, $M_t$-free graphs are polynomially $\eta$-bounded.
\end{conjecture}

This paper is organized as follows. In Section~\ref{sec:prel}, we introduce the required notation and terminology to be used in our proofs. Sections~\ref{sec:P51} and \ref{sec:P52} are devoted to the proof of Theorem~\ref{mainP5}. In Section~\ref{sec:stars}, we prove Theorems~\ref{mainKrs1} and \ref{mainFt1}. Finally, in Section~\ref{sec:poly}, we prove Theorems~\ref{properP51} and \ref{mainparallel1}.

\section{Notation and Terminology}\label{sec:prel}

Let $G = (V(G),E(G))$ be a graph. For a set $X \subseteq V(G)$ we denote by $G[X]$ the subgraph of $G$ induced by $X$. For $X \subseteq V(G)\cup E(G)$, $G \setminus X$ denotes the subgraph of $G$ obtained by removing $X$. Note that if $X\subseteq V(G)$, then $G \setminus X$ denotes the subgraph of $G$ induced by $V(G)\setminus X$.  In this paper, we use induced subgraphs and their vertex sets interchangeably.

For every $x\in G$, we denote by $N_G(x)$ the set of neighbours of $x$ in $G$, and write $N_G[x]=N_G(x)\cup \{x\}$ (we omit the subscript $G$ if there is no ambiguity). For an induced subgraph $H$ of $G$ (where $x$ does not necessarily belongs to $H$), we define $N_H(x)=N_G(x) \cap H$ and $N_H[x]=N_G[x]\cap H$. Also, for $X\subseteq G$, we denote by $N_G(X)$ the set of all vertices in $G\setminus X$ with at least one neighbour in $X$, and define $N_G[X]=N_G(X)\cup X$. 

Let $X,Y \subseteq G$ be disjoint. We say $X$ is \textit{complete} to $Y$ if all edges with an end in $X$ and an end in $Y$ are present in $G$, and $X$ is \emph{anticomplete}
to $Y$ if there are no edges between $X$ and $Y$. For $x\in V(G)\setminus Y$, we also say $x$ is \textit{complete} (\textit{anticomplete}) to $Y$ if $\{x\}$ is complete (anticomplete) to $Y$.

An {\em induced path in $G$} is an induced subgraph of $G$ that is a path. If $P$ is an induced path in $G$, we write $P = p_1 \dd \cdots \dd p_k$ to mean that $V(P) = \{p_1, \dots, p_k\}$ and $p_i$ is adjacent to $p_j$ if and only if $|i-j| = 1$. We call the vertices $p_1$ and $p_k$ the \emph{ends of $P$}, and say that $P$ is \emph{from $p_1$ to $p_k$}. The \emph{interior of $P$}, denoted by $P^*$, is the set $P \setminus \{p_1, p_k\}$. The \emph{length} of a path is its number of edges (so a path of length at most one has empty interior).

\section{Cradles in $P_5$-free graphs}\label{sec:P51}
We now plunge into proving Theorem~\ref{mainP5}. Roughly, the proof consists of the following two steps. First, we show that in every $P_5$-free graph $G$ of bounded clique number, one may always find a small hitting set for all maximum stable sets of $G$ which are ``restricted to'' a configuration in $G$ that we call a ``cradle''. The second step then is to show that in every $P_5$-free graph of bounded clique number, one may ``pack'' only a few cradles such that every maximum stable set of $G$ is restricted to at least one of the cradles from the packing. This, in particular, is done starting with a maximum clique $K$ in $G$ along with a linear order of the vertices in $K$, and then hitting maximum stable sets in $G$ intersecting $K$ and those not intersecting $K$ separately. The former is easily handled since $K$ has bounded cardinality. For maximum stable sets not intersecting $K$, we further group them based on the first, if any, vertex $v$ in the linear order of $K$ with a ``private'' neighbour in the stable set, and show that the stable sets in each group are restricted to a specific cradle corresponding to $v$.

In this section, we take the first step, and the second one is postponed to the next section. The main result is Theorem~\ref{hittingsoothed}, but we need some preparation before we state and prove it.

Let $G$ be a graph. By a \textit{cradle} in $G$ we mean a pair $(X,Z)$ of disjoint subsets of $V(G)$ such that either $|X|\leq 1$, or the following hold.
\begin{itemize}
    \item Every vertex in $N_Z(X)$ has a neighbour in $G\setminus N[X]$.
    \item For every two vertices $z,z'\in N_Z(X)$, there exists an induced path $P$ in $G$ from $z$ to $z'$ such that $P^*\subseteq G\setminus N[X]$.
\end{itemize} 
In particular, for every $X\subseteq V(G)$, $(X,\emptyset)$ is a cradle in $G$. The following lemma is immediate from the above definition:
\begin{lemma}\label{cradlematching}
Let $G$ be a $P_5$-free graph and let $(X,Z)$ be a cradle in $G$. Then the following hold.
\begin{itemize}
    \item For all distinct vertices $x, x' \in X$ and $z,z' \in Z$, we have $E(G[\{x, x', z, z'\}]) \neq \{xz, x'z'\}$.
    \item For every vertex $z\in Z$, there is at most one component of $G[X]$ in which $z$ has both a neighbour and a non-neighbour.
\end{itemize}
\end{lemma}
\begin{proof}
    Suppose the first bullet does not hold. Since $(X,Z)$ is a cradle, there exists an induced path $P$ of length at least two in $G$ from $z$ to $z'$ such that $P^*\subseteq G\setminus N[X]$. But now $x\dd z\dd P\dd z'\dd x'$ is an induced path in $G$ on at least five vertices, a contradiction.

    To see the second bullet, suppose for a contradiction that there exist distinct components $D,D'$ of $G[X]$ and induced paths $z\dd a\dd b$ and $z\dd a'\dd b'$ in $G$ such that $\{a,b\}\subseteq D$ and $\{a',b'\}\subseteq D$. But now $b\dd a\dd z\dd a'\dd b'$ is an induced $P_5$ in $G$, which is impossible. This proves Lemma~\ref{cradlematching}. 
\end{proof}

Given a cradle $(X,Z)$ in a graph $G$, by a \textit{rocker for} $(X,Z)$ we mean a pair $(\mathcal{I},\mathcal{J})$ of collections of components of $G[X]$ with the following specifications.
\begin{itemize}
    \item $\mathcal{I}$ ``minimally represents'' all components of $G[X]$ which are anticomplete to some vertex in $Z$. More precisely, for every vertex $z\in Z$ which is anticomplete to some component of $G[X]$, there exists $D\in \mathcal{I}$ such that $z$ is anticomplete to $D$, and subject to this property, $\mathcal{I}$ is minimal with respect to inclusion. 
    \item $\mathcal{J}$ is the collection of all components $D$ of $G[X]$ for which there exists a vertex $z \in Z$ such that $z$ has a neighbour in each component of $G[X]$, and $z$ has a non-neighbour in $D$. 
\end{itemize}
It follows that there exists a rocker for every cradle in $G$. Also, we deduce:
\begin{lemma}\label{rockersize}
    Let $G$ be a $P_5$-free graph. Let $(X,Z)$ be a cradle in $G$ and let $(\mathcal{I},\mathcal{J})$ be a rocker for $(X,Z)$. Then $|\mathcal{I}|,|\mathcal{J}|\leq \omega(G)$.
\end{lemma}
\begin{proof}
    First, we show that $|\mathcal{I}|\leq \omega(G)$. By the minimality of $\mathcal{I}$, for every $D \in \mathcal{I}$, there exists a vertex $z_D \in Z$ such that $z_D$ is anticomplete to $D$ but $z_D$ has a neighbour $x_{D,D'}$ in $D'$ for each $D' \in \mathcal{I} \setminus \{D\}$. It follows that the vertices in $Z'=\{z_D:D\in \mathcal{I}\}$ are pairwise distinct, and so $|\mathcal{I}|=|Z'|$. Also, for all distinct $D,D'\in \mathcal{I}$, from the first bullet of Lemma~\ref{cradlematching} applied to $(X,Z)$ and $x=x_{D,D'},x'=x_{D',D},z=z_{D}$ and $z'=z_{D'}$, we deduce that $z_D$ is adjancent to $z_{D'}$. Therefore, $Z'$ is a clique of $G$ and so $|\mathcal{I}|=|Z'|\leq \omega(G)$.
    
    The proof of  $|\mathcal{J}|\leq \omega(G)$ is similar. For every $D \in \mathcal{J}$, let $z'_D \in Z$ be a vertex with a neighbour in every component of $G[X]$ and with a non-neighbour $x_{D}$ in $D$. For distinct $D,D'\in \mathcal{J}$, it follows from the second bullet of Lemma~\ref{cradlematching} applied to $(X,Z)$ that $z'_D,z'_{D'}$ are distinct, $z'_D$ is adjacent to $x_{D'}$ and  $z'_{D'}$ is adjacent to $x_{D}$. In particular, assuming $Z''=\{z'_D : D \in \mathcal{J}\}$, we have $|\mathcal{J}|=|Z''|$. Now, by the first bullet of Lemma~\ref{cradlematching} applied to $(X,Z)$ and $x=x_D, x'=x_{D'}, z=z'_{D'}$ and $z'=z'_{D}$, $z'_D$ is adjacent to $z'_{D'}$. Therefore, $Z''$ is a clique of $G$ and so $|\mathcal{J}|=|Z''|\leq \omega(G)$. This completes the proof of Lemma~\ref{rockersize}.
\end{proof}

 Let $G$ be a graph and let $S\subseteq V(G)$. We say $S$ is \textit{$(X,Z)$-restricted} if $S\subseteq X\cup Z$ and $S\cap X\neq \emptyset$. Next we prove a lemma about $(X,Z)$-restricted maximum stable sets in $G$ (although it remain true for ``maximal stable sets'' as well).

\begin{lemma}\label{Dsoothed}
    Let $G$ be a $P_5$-free graph. Let $(X,Z)$ be a cradle in $G$ such that some vertex in $Z$ is not complete to $X$. Let $(\mathcal{I},\mathcal{J})$ be a rocker for $(X,Z)$. Then for every $(X,Z)$-restricted maximum stable set $S$ of $G$, we have $S\cap Q\neq \emptyset$ for some $Q\in \mathcal{I}\cup \mathcal{J}$.
\end{lemma}
\begin{proof}
Suppose that $S$ is an $(X,Z)$-restricted maximum stable set of $G$ such that $S\cap Q=\emptyset$ for all $Q\in \mathcal{I}$. Our goal is to show that $S\cap Q\neq \emptyset$ for some $Q\in \mathcal{Q}$. The assumption that some vertex in $Z$ is not complete to $X$, along with the second bullet of Lemma~\ref{cradlematching}, implies that $\mathcal{I}\cup \mathcal{J}\neq \emptyset$. Also, for every $Q\in \mathcal{I}$, since $S$ is a maximum stable set of $G$ which is disjoint from $Q$, it follows that there exists a vertex $z_Q\in S\setminus X=S\cap Z$ with a neighbour $x_Q\in Q$. In particular, we have $S\cap Z\neq \emptyset$. We claim that:

\sta{\label{seesall}There exists $z\in S\setminus X=S\cap Z$ with a neighbour in each component of $G[X]$.}

Suppose not. Let $z\in S\setminus X=S\cap Z$ be a vertex which has a neighbour in as many of the components of $G[X]$ as possible. Then $z$ is anticomplete to some component of $G[X]$. It follows from the choice of $\mathcal{I}$ that $z$ is anticomplete to some component $Q\in \mathcal{I}$. Note that since $z_{Q}\in S\cap Z$ is adjacent to $x_{Q}\in Q$, $z$ and $z_Q$ are distinct. Also, by the choice of $z$, there exists a component $Q'$ of $G[X]$ such that $z$ has a neighbour $x$ in $Q'$ and $z_Q$ is anticomplete to $Q'$. It follows that $Q$ and $Q'$ are distinct. But now, assuming $x'=x_Q$ and $z'=z_Q$, we have $E(G[\{x, x', z, z'\}])= \{xz, x'z'\}$, which violates the first bullet of Lemma~\ref{cradlematching}. This prove \eqref{seesall}.

\medskip

Let $z\in S\setminus X=S\cap Z$ be as promised by  \eqref{seesall}. Since $S \cap X\neq \emptyset$, it follows that there exists a vertex $z'\in S\cap Q$ for some component $Q$ of $G[X]$. But then $z$ has a non-neighbour (namely, $z'$) in $Q$, and so $Q\in \mathcal{J}$. This, along with $z'\in S\cap Q$, implies that $S \cap Q\neq \emptyset$ for $Q\in \mathcal{J}$, and completes the proof of Lemma~\ref{Dsoothed}.
\end{proof}

Let us now prove the main result of this section:

\begin{theorem}\label{hittingsoothed}
   For all integers $c,d,h\in \mathbb{N}$, there exists $\Gamma(c,d,h)\in \mathbb{N}$ with the following property. Let $G$ be a $P_5$-free graph with $\omega(G)\leq c$. Let $(X,Z)$ be a cradle in $G$ with $\omega(X)\leq d$ and $\eta(X')\leq h$ for every $X'\subseteq X$. Then there exists $W\subseteq G$ with $|W|\leq \Gamma(c,d,h)$ which intersects every $(X,Z)$-restricted maximum stable set of $G$.
\end{theorem}
\begin{proof}
We define $\Gamma(c,d,h)$ recursively, as follows. Let $c,h\in \mathbb{N}$ be arbitrary yet fixed. Let $\Gamma(c,1,h)=2c$, and for every $d\geq 2$, let
$$\Gamma(c,d,h)=2c(h + (c+1)^{c+1} \Gamma(c,d-1,h)).$$
 We prove by induction on $d$, with $c,h$ fixed, that $\Gamma(c,d,h)$ satisfies Theorem~\ref{hittingsoothed}.
 
 Let $\mathcal{S}$ be the sets of all $(X,Z)$-restricted maximum stable sets of $G$. Let $(\mathcal{I},\mathcal{J})$ be a rocker for $(X,Z)$. By Lemma~\ref{rockersize}, we have $|\mathcal{I}\cup\mathcal{J}|\leq 2c$. For each $Q\in \mathcal{I}\cup\mathcal{J}$, let $\mathcal{S}_Q$ be the set of all $(X,Z)$-restricted maximum stable sets of $G$ with $S\cap Q\neq \emptyset$. By Lemma~\ref{Dsoothed}, we have $\mathcal{S}=\bigcup_{Q\in \mathcal{I}\cup\mathcal{J}}\mathcal{S}_Q$.
 
 To launch the induction, consider the base case $d=1$. Then $X$ is a stable in $G$, and so every component of $G[X]$ is a singleton. But now $W=\bigcup_{Q\in \mathcal{I}\cup\mathcal{J}}Q$ is a set of at most $2c=\Gamma(c,1,h)$ vertices which intersects every $(X,Z)$-restricted maximum stable set  of $G$, as desired.
 
 From now on, assume that $d\geq 2$, and so $|X|\geq 2$.  
 Note that if every vertex in $Z$ is complete to $X$, then every set $S\in \mathcal{S}$ is contained in $X$. But then we are done as $\eta(X)\leq h\leq \Gamma(c,d,h)$. Thus, we may assume from now on that some vertex in $Z$ is not complete to $X$. Moreover, since $\mathcal{S}=\bigcup_{Q\in \mathcal{I}\cup\mathcal{J}}\mathcal{S}_Q$, in order to prove Theorem~\ref{hittingsoothed}, it suffices to show that for every $Q\in \mathcal{I}\cup\mathcal{J}$, there exists $W_Q\subseteq V(G)$ with $|W_Q|\leq h + (c+1)^{c+1} \Gamma(c,d-1,h)$ which intersects every set $S\in \mathcal{S}_Q$. Henceforth, let $Q\in \mathcal{I}\cup \mathcal{J}$ be fixed.

By the assumption, we have $\eta(Q)\leq h$, and so there exists a hitting set $W_0\subseteq Q$ for $G[Q]$ with $|W_0|\leq h$. Let us call a set $S\in \mathcal{S}_Q$ is \emph{mischievous} if $S \cap W_0 = \emptyset$. We deduce:

\sta{\label{obtainnbr} For every mischievous set $S\in \mathcal{S}_Q$, we have $S\cap N_Z(Q)\neq \emptyset$.}

Note that since $S$ is mischievous, it follows that $S \cap Q$ is not a maximum stable set of $G[Q]$. Therefore, since $S$ is an $(X,Z)$-restricted maximum stable set of $G$, we have $S\cap N_{X\cup Z}(Q)\neq \emptyset$. But $Q$ is component of $G[X]$, and so $N_{X}(Q)=\emptyset$. This yields $S\cap N_Z(Q)\neq \emptyset$, and so proves \eqref{obtainnbr}.\medskip

Next, let $Y\subseteq N_Z(Q)$ be chosen such that $Y$ ``minimally represents'' all vertices in $N_Z(Q)$ which are not complete to $Q$. More precisely, let $Y$ be a minimal subset (with respect to inclusion) of $N_Z(Q)$ such that every vertex in $Q$ with a non-neighbour in $N_Z(Q)$ has a non-neighbour in $Y$. It turns out that $Y$ is small:

\sta{\label{Ysmall}We have $|Y|<(c+1)^{c+1}$.}

Suppose for a contradiction that $|Y| \geq (c+1)^{c+1}$. Since $\omega(G)\leq c$, it follows from Theorem~\ref{classicalramsey} that $Y$ contains a stable set $Y'$ of cardinality $c+1$. Also, note that by the minimality of $Y$, for every $y \in Y$, there is a vertex $q_y \in Q$ such that $y$ is the only non-neighbour of $q_y$ in $Y$. Therefore, again since $\omega(G)\leq c$, there exist two distinct vertices $y, y' \in Y'$ such that $q_y$ and $q_{y'}$ are non-adjacent. But now $x=q_{y}, x'=q_{y'}, z=y'$ and $z'=y$ violate the first bullet of Lemma~\ref{cradlematching} applied to $(X,Z)$. This proves \eqref{Ysmall}.\medskip 

Now, for every $y\in Y$, let
$$X_y= Q\setminus N_Q(y),$$
$$Z_y = (X\setminus X_y)\cup Z.$$
It follows that:

\sta{\label{obtaincradle} For every $y\in Y$, $(X_y,Z_y)$ is a cradle in $G$.}

From the definition of $X_y$ and $Z_y$, it follows immediately that $X_y\cap Z_y=\emptyset$, and $N_{Z_y}(X_y)\subseteq N_Q(y)\cup  Z$. Note that every vertex in $N_Q(y)$ is adjacent to $y\in G\setminus N[X_y]$. Also, since $(X,Z)$ is a cradle in $G$, every vertex in $N_Z(X)$ has a neighbour in $G\setminus N[X]\subseteq G\setminus N[X_y]$. We conclude that every vertex in $N_{Z_y}[X_y]\subseteq N_Q(y)\cup  N_Z(X)$ has a neighbour in $G\setminus N[X_y]$. It remains to show that for every two distinct vertices $z,z'\in N_{Z_y}[X_y]\subseteq N_Q(y)\cup  N_Z(X)$, there exists an induced path $P$ in $G$ from $z$ to $z'$ such that $P^*\subseteq G\setminus N[X_y]$. If $z,z'\in N_Q(y)$, then there exists an induced path $P$ in $G$ from $z$ to $z'$ with $P^*\subseteq \{y\}\subseteq G\setminus N[X_y]$. Also, if $z,z'\in N_Z(X)$, then since $(X,Z)$ is a cradle, it follows that there exists an induced path $P$ in $G$ from $z$ to $z'$ with $P^*\subseteq G\setminus N[X]\subseteq G\setminus N[X_y]$. So we may assume that $z\in N_Q(y)$ and $z'\in N_Z(X)$. Since $(X,Z)$ is a cradle in $G$ and $y,z'\in N_Z(X)$, there exists an induced path $P'$ in $G$ from $y$ to $z'$ with $P'^*\subseteq G\setminus N[X]$. Moreover, $y,z$ are adjacent in $G$. Now $P=z\dd y\dd P'\dd z'$ is an induced path in $G$ from $z$ to $z'$ with $P^*\subseteq (G\setminus N[X])\cup \{y\}\subseteq G\setminus N[X_y]$. This proves \eqref{obtaincradle}.

\medskip

We also need to investigate how mischievous sets interact with the cradles found in \eqref{obtaincradle}. For every $y\in Y$, let $\mathcal{S}^y$ be the set of all mischievous sets $S\in \mathcal{S}_Q$ which are $(X_y,Z_y)$-restricted.

\sta{\label{obtainfriend} For every mischievous set $S\in \mathcal{S}_Q$, there exists $y \in Y$ such that $S\in \mathcal{S}^y$.}

 Since $S$ is $(X,Z)$-restricted, it follows that $S\subseteq X\cup Z=X_y\cup Z_y$ for all $y\in Y$. We now need to show that $S\cap X_y\neq \emptyset$ for some $y\in Y$. Since $S\in \mathcal{S}_Q$, we may choose $q\in S\cap Q$. Also, since $S$ is mischievous, by \eqref{obtainnbr}, there exists a vertex $z\in S\cap N_Z(Q)$. It follows that $q\in S\cap Q$ has a non-neighbour in $N_Z(Q)$ (namely $z$), and so by the definition of $Y$, $q$ has a non-neighbour $y\in Y$. As a result, we have $q\in (S\cap Q) \setminus N_Q(y)=S\cap X_y$, and so $S\cap X_y\neq \emptyset$. This proves \eqref{obtainfriend}.

\medskip

Recall that the proof is by induction on $d$. In order to apply the induction hypothesis, we need the following two statements.

\sta{\label{prepforih} Let $y\in Y$. Then for every $q \in N_Q(y)$ and every component $D$ of $X_y$, $q$ is either complete or anticomplete to $D$.} 

For otherwise, there is an induced path $q\dd a\dd b$ in $G$ with $a, b \in D$. Since $y \in Z$ is adjacent to $q\in Q\subseteq X$, we have $y\in N_Z[X]$. Thus, since $(X,Z)$ is a cradle in $G$, it follows that $y$ has a neighbour in $u\in G\setminus N[X]$. But now $u\dd y\dd q\dd a\dd b$ is an induced $P_5$ in $G$, a contradiction. This proves \eqref{prepforih}.

\sta{\label{almostih} For every $y\in Y$, we have $\omega(X_y) \leq d-1$.}

Suppose not. Let $K$ be a clique of cardinality $d$ in $X_y$. Let $D$ be the component of $X_y$ with $K\subseteq D$. Since $Q$ is connected, there exists a vertex $q\in Q\setminus X_y=N_Q(y)$ with a neighbour in $D$. By \eqref{prepforih}, $q$ is complete to $D$. But then $K\cup \{q\}$ is a clique in $Q\subseteq X$ of cardinality $d+1$, which violates the assumption $\omega(X)\leq d$. This proves \eqref{almostih}.

\medskip

\sta{\label{ih} For every $y\in Y$, there exists $W^y\subseteq G$ with $|W^y|\leq \Gamma(c,d-1,h)$ which intersects every set $S\in \mathcal{S}^y$.}

By \eqref{obtaincradle}, \eqref{almostih}, and the fact that $X_y\subseteq X$, $(X_y,Z_y)$ is a cradle in $G$ with $\omega(X_y)\leq d-1$ and $\eta(X_y)\leq h$. Therefore, by the induction hypothesis applied to $(X_y,Z_y)$, there exists $W^y\subseteq G$ with $|W^y|\leq \Gamma(c,d-1,h)$ which intersects every mischievous set $S\in \mathcal{S}_Q$ which is $(X_y,Z_y)$-restricted. In other words, $W^y$ intersects every set $S\in \mathcal{S}^y$. This proves \eqref{ih}.

\medskip

The proof is almost concluded. Let $W_Q = W_0 \cup (\bigcup_{y \in Y} W^y)$, where $W^y$ is as in \eqref{ih}. Then by \eqref{Ysmall}, we have $|W| \leq h + (c+1)^{c+1} \Gamma(c,d-1,h)=\Gamma(c,d,h)$. Let $S\in \mathcal{S}_Q$ be arbitrary. We need to show that $S\cap W\neq \emptyset$. If $S\cap W_0\neq \emptyset$, then we have $S\cap W\neq \emptyset$ as $W_0\subseteq W$. So we may assume that $S \cap W_0=\emptyset$, that is, $S$ is mischievous. By \eqref{obtainfriend}, there exists $y\in Y$ such that $S\in \mathcal{S}^y$. By \eqref{ih}, we have $S\cap W^y\neq \emptyset$, which along with the fact that $W^y\subseteq W$, implies that $S\cap W\neq \emptyset$. This completes the proof of Theorem~\ref{hittingsoothed}. 
\end{proof}

\section{Cradle packing and proof of Theorem~\ref{mainP5}}\label{sec:P52}
Here we complete the proof of Theorem~\ref{mainP5}. First, note that for every maximum stable set in a graph $G$ and every connected component $G'$ of $G$, $S\cap G'$ is a maximum stable set of $G'$. This immediately yields the following.
\begin{lemma} \label{lem:conn}
For every graph $G$ and every connected component $G'$, we have $\eta(G) \leq \eta(G')$. 
\end{lemma}

Let $G$ be a graph. For a cradle $\xi=(X,Z)$ in $G$, we write $X_{\xi}$ for $X$ and $Z_{\xi}$ for $Z$.  A \textit{cradle packing in $G$} is a collection $\Xi$  of distinct cradles in $G$ such that every maximum stable set of $G$ is $\xi$-restricted for some $\xi\in \Xi$. We also define $\eta(\Xi)=\max\{\eta(X'):X'\subseteq X_{\xi}, \ \xi \in \Xi\}$. Then we have: 

\begin{lemma}\label{cradlepacking}
For all $c,h,t\in \mathbb{N}$, there exists $\Phi(c,h,t)\in \mathbb{N}$ with the following property. Let $G$ be a $P_5$-free graph with $\omega(G)\leq c$. Assume that there exists a cradle packing $\Xi$ in $G$ with $|\Xi|\leq t$ and $\eta(\Xi)\leq h$. Then we have $\eta(G)\leq \Phi(c,h,t)$.
\end{lemma}
\begin{proof}
Let $\Gamma(\cdot,\cdot,\cdot)$ be as in Theorem~\ref{hittingsoothed} and let  $\Phi(c,h,t)=t\Gamma(c,c,h)$. It follows from Theorem~\ref{hittingsoothed} that, for every $\xi \in \Xi$, there exists $W_{\xi}\subseteq G$ with $|W_{\xi}|\leq \Gamma(c,c,h)$ which intersects every $\xi$-restricted maximum stable set of $G$. Now, let $W=\bigcup_{\xi\in \Xi} W_{\xi}$. Then since $|\Xi|\leq t$, we have $|W|\leq t\Gamma(c,c,h)=\Phi(c,h,t)$. It remains to show that $W$ is a hitting set for $G$. Let $S$ be a maximum stable set of $G$. Since $\Xi$ is a cradle packing in $G$, it follows that $S$ is $\xi$-restricted for some $\xi\in \Xi$. But then we have $S\cap W_{\xi}\neq \emptyset$, which along with the fact that $W_{\xi}\subseteq W$, implies that $S\cap W\neq \emptyset$.  This completes the proof of Theorem~\ref{cradlepacking}. 
\end{proof}

Here comes the proof of Theorem~\ref{mainP5}, which we restate:

\begin{theorem}\label{mainP5proof}
For every $c\in \mathbb{N}$, there exists $\Psi(c)\in \mathbb{N}$ such that every $P_5$-free graph with $\omega(G)\leq c$ satisfies $\eta(G)\leq \Psi(c)$.
\end{theorem}

\begin{proof}
Let $\Psi(1)=1$, and for every $c\geq 2$, let
$\Psi(c)=\Phi(c,\Psi(c-1),2c+1)$,
where $\Phi(\cdot,\cdot,\cdot)$ is as in Lemma~\ref{cradlepacking}. We prove by induction on $c$ that  every $P_5$-free graph with $\omega(G)\leq c$ satisfies $\eta(G)\leq \Psi(c)$. The result is easily seen to hold for $\omega(G)\leq 1$. Thus, writing $\omega(G)=k$, we may assume that $c\geq k\geq 2$. Also, by Lemma \ref{lem:conn}, we may assume that $G$ is connected.

Let $K=\{v_1,\ldots, v_k\}$ be a maximum clique in $G$. For every $i\in \{1,\ldots, k\}$, let us define $\zeta_i=(\{v_i\},G\setminus \{v_i\})$; then $\zeta_i$ is a cradle in $G$ because the first set in the pair $\zeta_i$ is a singleton (roughly, these cradles will be used to hit the maximum stable sets of $G$ which intersect $K$).

We now define another set of cradles in $G$. Let $X_0$ be the set of all vertices in $G\setminus K$ which are anticomplete to $K$. For each $i\in \{1,\ldots, k\}$, let $X_i$ be the set of all vertices in $G\setminus K$ which are adjacent to $v_i$ and non-adjacent to $v_j$ for all $j\in \{1,\ldots, k\}\setminus \{1,\ldots, i\}$. It follows that $X_0,X_1,\ldots,X_k$ is a partition of $G\setminus K$. For every $i\in \{0,1,\ldots, k-1\}$, let $Z_i=\bigcup_{j=i+1}^k X_i$, and let $Z_k=\emptyset$. Define $\xi_i=(X_i,Z_i)$ for each $i\in \{0,1,\ldots, k\}$. We claim that:

\sta{\label{st:cradles} For all $i \in \{0,1,\ldots, k\}$, $\xi_i$ is a cradle in $G$.}

The assertion is trivial for $i=k$ as $Z_k=\emptyset$. Also, for each $i \in \{0,1,\ldots, k-1\}$, note that $K_i=\{v_{i+1},\ldots, v_k\}$ is a clique in $G$ with $K_i\subseteq G\setminus N[X_i]$, and every vertex in $Z_i$ has a neighbour in $K_i$. This proves \eqref{st:cradles}. 

\medskip

From \eqref{st:cradles} and the definition of $X_0,X_1,\ldots, X_k$, it is immediately seen that for every non-empty subset $S\subseteq V(G)$, either $S\cap K\neq \emptyset$, or $S$ is $\xi_i$-restricted where $i=\min\{j:i\in \{0,\ldots,k\},S\cap X_j\neq \emptyset\}$. As a result, every non-empty subset of $V(G)$ is either $\zeta_i$-restricted for some $i\in \{1,\ldots,k\}$ or $\xi_i$-restricted for some $i\in \{0,1,\ldots,k\}$. Consequently, defining 
$$\Xi=\{\zeta_i:i\in \{1,\ldots,k\}\}\cup \{\xi_i:i\in \{0,1,\ldots,k\}\},$$
it follows that $\Xi$ is a cradle packing in $G$ with $|\Xi|\leq 2k+1\leq 2c+1$. Moreover, we deduce:

\sta{\label{st:clique} $\eta(\Xi)\leq \Phi(c-1)$.} 

From the definition of $\zeta_1,\ldots, \zeta_k$, we only need to show that for every $i\in \{0,1,\ldots,k\}$, and every connected component $D$ of $G[X_i]$, we have  $\eta(D)\leq \Phi(c-1)$. To that end, using the induction hypothesis, it suffices to show that $\omega(D)\leq c-1$. Suppose for a contradiction that there exists a clique $K'$ in $G$ of cardinality $c$ with $K'\subseteq D$. Then we have $i=0$, as otherwise $K'\cup \{v_i\}$ would be a clique of cardinality $c+1>\omega(G)$ in $G$. Since $G$ is connected and $K'\subseteq X_0$ is anticomplete to $K$, there exist an induced path $P$ in $G$ with (not necessarily distinct) ends $x$ and $x'$ such that $x$ has a neighbour $y$ in $K$, $x'$ has a neighbour $y'$ in $K'$, $P\setminus \{x\}$ is anticomplete to $K'$  and $P\setminus \{x'\}$ is anticomplete to $K$. Since both $K$ and $K'$ are maximum cliques of $G$, it follows that $x$ has a non-neighbour $z\in K$ and $x'$ has a non-neighbour $z' \in K'$. But now $z'\dd y'\dd x'\dd P\dd x\dd y\dd z$ is an induced path in $G$ on at least five vertices, a contradiction. This proves \eqref{st:clique}. 
\medskip

To sum up, for the $P_5$-free graph $G$ with $\omega(G)\leq c$, we proved that $\Xi$ is a cradle packing in $G$ for which we have $|\Xi|\leq 2c+1$, and, by \eqref{st:clique}, $\eta(\Xi)\leq \Psi(c-1)$. But then Lemma~\ref{cradlepacking} implies that $\eta(G)\leq \Phi(c,\Psi(c-1),2c+1)=\Phi(c)$.
This completes the proof of Theorem~\ref{mainP5proof}. 
 \end{proof}

\section{More than stars}\label{sec:stars}
In this section, we prove Theorems~\ref{mainKrs1} and \ref{mainFt1}, beginning with following lemma.

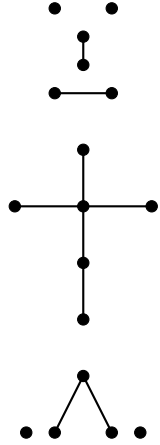
\begin{figure}
\centering

\begin{tikzpicture}[scale=1.5,auto=left]
\tikzstyle{every node}=[inner sep=1.5pt, fill=black,circle,draw]  
\centering
\node at (0.5,0.5) {};
\node at (-0.5,0.5) {};
\node at (-0.25,0.5) {};
\node at (0.25,0.5) {};
\node at (0,1) {};
\node at (0,1.5) {};
\node at (0,2) {};
\node at (0,2.5) {};
\node at (0,3) {};
\node at (0.6,2.5) {};
\node at (-0.6,2.5) {};
\node at (0.25,3.5) {};
\node at (-0.25,3.5) {};
\node at (0,3.75) {};
\node at (0,4) {};
\node at (0.25,4.25) {};
\node at (-0.25,4.25) {};

\draw [thick] (0,4) -- (0,3.75);
\draw [thick] (0.25,3.5) -- (-0.25,3.5);

\draw [thick] (0,1.5) -- (0,2);
\draw [thick] (0,2) -- (0,2.5);
\draw [thick] (0,2.5) -- (0,3);
\draw [thick] (0,2.5) -- (0.6,2.5);
\draw [thick] (0,2.5) -- (-0.6,2.5);

\draw [thick] (0,1) -- (-0.25,0.5);
\draw [thick] (0,1) -- (0.25,0.5);

\end{tikzpicture}

\caption{Graphs $L_2$ (top), $F_3$ (middle) and $S_{(2,2)}$ (bottom).}
\label{fig:graphs}
\end{figure}

\begin{lemma}\label{hitmanytimes}
    Let $G$ be a graph and let $k\in \mathbb{N}$. For every $i\in \{1,\ldots, k\}$, let $d_i,h_i$ be two positive integers, and let $A_i,A'_i$ be two subsets of $V(G)$, such that the following hold.
    \begin{itemize}
    \item We have $A_i\cup A'_i=V(G)$ for all $i\in \{1,\ldots, k\}$.
        \item For each $i\in \{1,\ldots, k\}$ and every $X\subseteq A_i$, we have $\eta(X)\leq h_i$.
        \item For every maximum stable set $S$ of $G$, we have $|S\cap A'_i|< d_i$ for some $i\in \{1,\ldots, k\}$.
    \end{itemize}
    Then $\eta(G)\leq \Sigma_{i=1}^kd_ih_i$.
\end{lemma}
\begin{proof}
From Lemma~\ref{etatimes} applied to $G[A_i]$, $d_i$ and $h_i$, we deduce that:

\sta{\label{getB} For each $i\in \{1,\ldots,k\}$, there exists $D_i\subseteq A_i$ with $|D_i|\leq d_ih_i$ and $\alpha(A_i\setminus D_i)\leq \alpha(A_i)-d_i\leq \alpha(G)-d_i$.}

   Now, let $D=D_1\cup \cdots \cup D_k$. We claim that $D$ is a hitting set for $G$. Suppose for a contradiction that there exists a maximum stable set $S$ of $G$ with $S\cap D=\emptyset$. By the third bullet of Lemma~\ref{hitmanytimes}, we have $|S\cap A'_i|<d_i$ for some $i\in \{1,\ldots,k\}$. On the other hand, we have $S\cap A_i\subseteq A_i\setminus D\subseteq A_i\setminus D_i$. Thus, from \eqref{getB}, it follows that $|S\cap A_i|\leq \alpha(A_i\setminus D_i)\leq \alpha(G)-d_i$. But then from the first bullet of Lemma~\ref{hitmanytimes}, we deduce that $|S|\leq |S\cap A_i|+|S\cap A_i'|<\alpha(G)$, a contradiction. The claim follows. Moreover, by \eqref{getB}, we have $|D|\leq |D_1|\cdots +|D_k|\leq d_1h_1+\cdots +d_kh_k$. This proves Lemma~\ref{hitmanytimes}. 
\end{proof}
Let us now prove Theorem~\ref{mainKrs1}, restated as follows. For all integers $s\geq 1$ and $t\geq 0$, let $S_{(s,t)}$ denote the unique graph (up to isomorphism) on $1+s+t$ vertices with a vertex of degree $s$, $s$ vertices of degree one and $t$ vertices of degree zero. Note that $S_{(s,0)}$ is isomorphic to the star $K_{1,s}$ (see also Figure~\ref{fig:graphs} for a depiction of $S_{(2,2)}$).
\begin{theorem}\label{mainKrs}
For all integers $c,s\geq 1$ and $t\geq 0$, there exists $\Psi(c,s,t)\in \mathbb{N}$ such that every $S_{(s,t)}$-free graph $G$ with $\omega(G)\leq c$ satisfies $\eta(G)\leq \Psi(c,s,t)$. Moreover, for all $c,s\in \mathbb{N}$, we have $\Psi(c,s,1)\leq c^{2s+1}$.
\end{theorem}
\begin{proof}
 Let $s\geq 1$ be fixed. For every $c\geq 1$, let $\Psi(c,s,0)=c^s$, and for every $t\geq 0$, let $\Psi(1,s,t)=1$. For every $c>1$ and $t>0$, let
$$\Psi(c,s,t)=t\Psi(c-1,s,t)+(s+1)\Psi(c,s,t-1).$$
We prove by induction on $c+t$ that every $S_{(s,t)}$-free graph $G$ with $\omega(G)\leq c$ satisfies $\eta(G)\leq \Psi(c,s,t)$. The result is trivial for $c=1$, and follows from Theorem~\ref{bndddegree} for $t=0$; thus, we may assume that $c>1$ and $s>0$. Let $x\in V(G)$ be arbitrary. Let us define $A_1=N(x)$, $A_2=G\setminus N[x]$, $A'_1=G\setminus N(x)$, $A'_2=N[x]$. We deduce:
  
  \sta{\label{st:lemmaconditions1}The following hold.
    \begin{itemize}
    \item We have $A_1\cup A'_1=A_2\cup A'_2=V(G)$.
        \item For every $X\subseteq A_1$, we have $\eta(X)\leq \Psi(c-1,s,t)$, and for every $X\subseteq A_2$, we have $\eta(X)\leq \Psi(c,s,t-1)$.
        \item For every maximum stable set $S$ of $G$, we have either $|S\cap A'_1|< t$ or $|S\cap A'_2|< s+1$.
    \end{itemize}}

    The first bullet is immediate. The second bullet follows directly from the induction hypothesis and the fact that $\omega(A_1)<\omega(G)\leq c$ and $G[A_2]$ is $S_{(s,t-1)}$-free (as otherwise $G[A_2\cup \{x\}]$ contains an induced $S_{(s,t)}$). To see the third bullet, suppose for a contradiction that there exists a maximum stable set $S$ of $G$ with $|S\cap A_1'|\geq t$ and $|S\cap A_2'|\geq s+1\geq 2$. It follows that $x\notin S$, and so  $|S\cap (G\setminus N[x])|\geq t$ and $|S\cap N(x)|\geq s+1$. Consequently, we can choose $A_1''\subseteq S\cap (G\setminus N[x])$ and $A_2''\subseteq S\cap N(x)$ with $|A_1''|=t$ and $|A_2''|=s$. But now $G[A_1''\cup A_2''\cup \{x\}]$ is isomorphic to $S_{(s,t)}$, a contradiction. This proves \eqref{st:lemmaconditions1}.

    \medskip

Now, \eqref{st:lemmaconditions1} allows an application of Lemma~\ref{hitmanytimes} to $d_1=t,d_2=s+1$, $h_1=\Psi(c-1,s,t), h_2=\Psi(c,s,t-1)$, $G,(A_1,A_1')$ and $(A_2,A_2')$. But then we have 
$$\eta(G)\leq t\Psi(c-1,s,t)+(s+1)\Psi(c,s,t-1)=\Psi(c,s,t).$$
Finally, let us show that $\Psi(c,s,1)\leq c^{2s+1}$ for all $c,s\in \mathbb{N}$. The proof is by induction on $c$ for fixed $s$. Note that $\Psi(1,s,1)=1$. For $c\geq 2$, we have $s+1\leq 2^s\leq c^s$. Combining this with the induction hypothesis and the definition of $\Psi(\cdot,\cdot,\cdot)$ yields:
\[\Phi(c,s,1)=\Psi(c-1,s,1)+(s+1)\Psi(c,s,0)\leq (c-1)^{2s+1}+c^{2s}<(c-1)c^{2s}+c^{2s}=c^{2s+1}.\]
This completes the proof of Theorem~\ref{mainKrs}.
\end{proof}
Next comes an easy lemma; we omit the proof (note also that this extends Theorem~\ref{bndddegree0}).
\begin{lemma}\label{cutsetlemma}
    Let $G$ be a graph, let $C\subseteq G$ and let $G'$ be a component of $G\setminus C$. Then $\eta(G)\leq |C|+\eta(G')$.
\end{lemma}
 We can now restate and prove Theorem~\ref{mainFt1}. For every $t\in \mathbb{N}$, let $F_t$ be the graph obtained from the star $K_{1,t+1}$ by subdividing an edge, exactly once. Equivalently, $F_t$ is the unique graph (up to isomorphism) on $t+3$ vertices with $t+1$ vertices of degree one, one vertex of degree two and one vertex of degree $t+1$ (see Figure~\ref{fig:graphs} for a depiction of $F_{3}$).
\begin{theorem}\label{mainFt}
For all $c,t\in \mathbb{N}$, there exists $\Psi(c,t)\in \mathbb{N}$ such that every $F_{t}$-free graph $G$ with $\omega(G)\leq c$ satisfies $\eta(G)\leq \Psi(c,t)$.
\end{theorem}
\begin{proof}
 Let $t\geq 1$ be fixed. Let $\Psi(1,t)=1$ and for every integer $c>1$, let
$$\Psi(c,t)=3(c+1)^{2t+1}+2t\Psi(c-1,t).$$
   We prove by induction on $c$ that the above value of $\Psi(c,t)$ satisfies Theorem~\ref{mainFt}. The result is trivial for $c=1$; thus we may assume that $c>1$.

Let $x_1\in V(G)$ be a vertex with at least one neighbour in $G$. Let $N=N_G(x_1)$ and $M=G\setminus N_G[x_1]$. Let $x_2\in N$ be chosen with $|N_M(x_2)|$ as large as possible. Let $P=N_M(x_2)$ and $Q=M\setminus P$. We begin with the following two statements, which have similar proofs.

   \sta{\label{st:N(Q)}We have $|N_Q(P)|< (c+1)^{2t+1}$.}
   Suppose not. Let $P'\subseteq P$ be minimal such that every vertex in $N_Q(P)$ has a neighbour in $P'$. Assume that $|P'|< (c+1)^{t+1}$. Then there exists a vertex $p\in P'$ such that $N_Q(p)> (c+1)^{t}$. Since $\omega(G)\leq c$, it follows from Theorem~\ref{classicalramsey} that there exists a stable set $A\subseteq N_Q(p)$ with $|A|=t$. But now $G[A\cup \{p,x_2\}]$ is isomorphic to $F_t$, a contradiction. We deduce that $|P'|>(c+1)^{t+1}$. Again, since $\omega(G)\leq c$, by Theorem~\ref{classicalramsey}, there exists a stable set $A'\subseteq P'$ with $|A'|=t+1$. Let $p'\in A'$ be arbitrary. By the minimality of $P'$, there exists a vertex $q\in N_Q(P)$ such that $N_{A'}(q)=\{p'\}$. But now $G[A'\cup \{q,x_2\}]$ is isomorphic to $F_t$, a contradiction. This proves \eqref{st:N(Q)}.

    \sta{\label{st:N(N)}We have $|N_Q(N)|< 2(c+1)^{2t+1}$.}
  Suppose not. Let $Q'=N_Q(N)\setminus N_Q(P)$. By \eqref{st:N(Q)}, we have $|Q'|>(c+1)^{2t+1}$. Let $N'\subseteq N$ be minimal such that every vertex in $Q'$ has a neighbour in $N'$. Assume that $|N'|< (c+1)^{t+1}$. Then there exists a vertex $y\in N'$ such that $N_{Q'}(y)> (c+1)^{t}$. Since $\omega(G)\leq c$, it follows from Theorem~\ref{classicalramsey} that there exists a stable set $A\subseteq N_{Q'}(y)$ with $|A|=t$. By the choice of $x_2$, there exists a vertex $p\in P$ which is not adjacent to $y$. Also, by the definition of $Q'$, $p$ is anticomplete to $Q'$. Now if $x_2,y$ are adjacent, then $G[A\cup \{p,x_2,y\}]$ is isomorphic to $F_t$, and if $x,y$ are not adjacent, then $G[A\cup \{x_1,x_2,y\}]$ is isomorphic to $F_t$, a contradiction. It follows that $|N'|>(c+1)^{t+1}$. Again, since $\omega(G)\leq c$, by Theorem~\ref{classicalramsey}, that there exists a stable set $A'\subseteq N'$ with $|A'|=t+1$. Let $y'\in A'$ be arbitrary. By the minimality of $N'$, there exist a vertex $q'\in Q'$ such that $N_{N'}(q')=\{y'\}$. But now $G[A'\cup \{q',x_1\}]$ isomorphic to $F_t$, which is impossible. This proves \eqref{st:N(Q)}.
  \medskip

  Let us define $A_1=N_G(x_1)=N$, $A_2=N_G(x_2)$, $A'_1=N_G(x_2)\setminus N_G(x_1)$, $A'_2=N_G(x_1)\setminus N_G(x_2)$ and $G'=A_1\cup A'_1=A_2\cup A_2'$. We deduce:
  
  \sta{\label{st:lemmaconditions2}The following hold.
    \begin{itemize}
        \item For each $i\in \{1,2\}$ and every $X\subseteq A_i$, we have $\eta(X)\leq \Psi(c-1,t)$.
        \item For every maximum stable set $S$ of $G$, we have either $|S\cap A'_1|< t$ or $|S\cap A'_2|< t$.
    \end{itemize}}

    The first bullet follows directly from the induction hypothesis and the fact that $\omega(A_1),\omega(A_2)<\omega(G)\leq c$. To see the second bullet, suppose for a contradiction that there exists a maximum stable set $S$ of $G$ with $|S\cap A_1'|,|S\cap A_2'|\geq t\geq 1$; consequently, we can choose $A_1''\subseteq S\cap A_1'$ and $A_2''\subseteq S\cap A_2'$ with $|A_1''|=t$ and $|A_2''|=1$. In addition, since $x_1$ is complete to $A_1'$ and $x_2$ is complete to $A_2'$ it follows that $x_1,x_2\notin S$. But now $G[A_1''\cup A_2''\cup \{x_1,x_2\}]$ is isomorphic to $F_t$, a contradiction. This proves \eqref{st:lemmaconditions2}.

    \medskip

To finish the proof, note that in view of \eqref{st:lemmaconditions2}, we can apply Lemma~\ref{hitmanytimes} to $d_1=d_2=t$, $h_1=h_2=\Psi(c-1,t)$, $G' (A_1,A_1')$ and $(A_2,A_2')$, which yields $\eta(G')\leq 2t\Psi(c-1,t)$. Furthermore, assuming $C=N_Q(P)\cup N_Q(N)$, it follows that $G'$ is a component of $G\setminus C$, and from \eqref{st:N(Q)} and \eqref{st:N(N)}, it is immediate that $|C|\leq 3(c+1)^{2t+1}$. Hence, by Lemma~\ref{cutsetlemma}, we have 
$$\eta(G)\leq |C|+\eta(G')\leq 3(c+1)^{2t+1}+2t\Psi(c-1,t)=\Psi(c,t).$$
This completes the proof of Theorem~\ref{mainFt}.
\end{proof}

\section{Proper induced subgraphs of $P_5$}\label{sec:poly}
In this last section, we prove Theorems~\ref{properP51} and \ref{mainparallel1} in reverse order. We need the following result of Wagon \cite{wagon}; recall that $M_t$ denotes the unique (up to isomorphism) $1$-regular graph on $2t$ vertices.

\begin{theorem}[Wagon \cite{wagon}]\label{wagon}
    For every $t\in \mathbb{N}$, every $M_t$-free graph $G$ satisfies $\chi(G)\leq \omega(G)^{2t-2}$.
\end{theorem}

We also need the following.

\begin{lemma}\label{incex}
    Let $s\geq 2$ be an integer and $X_1,\ldots, X_s$ be sets. Then we have $$\left|\bigcup_{i=1}^{s}X_i\right|\geq \sum_{i=1}^{s}\left|X_i\right|-\sum_{1\leq i<i'\leq s}|X_i\cap X_{i'}|.$$
\end{lemma}
\begin{proof}
    We induct on $s$. The case $s=2$ follows from the the inclusion-exclusion formula. For $s\geq 3$, by the inclusion-exclusion formula and the union bound, we have
    $$\left|\bigcup_{i=1}^{s}X_i\right|=\left|\bigcup_{i=1}^{s-1}X_i\right|+|X_s|-\left|\bigcup_{i=1}^{s-1}(X_i\cap X_s)\right|\geq \left|\bigcup_{i=1}^{s-1}X_i\right|+|X_s|-\sum_{i=1}^{s-1}\left|X_i\cap X_s\right|.$$
    Now the result is immediate from the induction hypothesis.
\end{proof}
Next we give a proof of Theorem~\ref{mainparallel1}, restated below. For every $t\in \mathbb{N}$, let $L_t$ be the unique graph (up to isomorphism) on $t+4$ vertices with $t$ vertices of degree zero and four vertices of degree one (see Figure~\ref{fig:graphs} for a depiction of $L_2$). 

\begin{theorem}\label{mainparallel}
For all $c,t\in \mathbb{N}$, every $L_{t}$-free graph $G$ with $\omega(G)\leq c$ satisfies $\eta(G)\leq c^{7t+7}$.
\end{theorem}
\begin{proof}
    The assertion is trivial for $c=1$; assume that $c>1$. Let $G$ be an $L_{t}$-free graph with $\omega(G)\leq c$. Since $M_{t+2}$ contains $L_t$ as an induced subgraph, by Theorem~\ref{wagon}, we have $\chi(G)\leq c^{2t+2}$, that is, $V(G)$ can be partitioned into $k$ stable sets $C_1,\ldots, C_k$ where $k\leq c^{2t+2}$. For each $i\in \{1,\ldots, k\}$, let $Z_i$ be the set of all vertices in $z\in V(G)\setminus C_i$ with $k|C_i\setminus N(z)|\geq |C_i|$.
    
 \sta{\label{st:smallD} Let $i\in \{1,\ldots, k\}$ such that $|C_i|\geq 4k^2t$. Then there exists a non-empty subset $D_i$ of  $C_i$ with $|D_i|<2k^2$ such that every vertex in $Z_i$ has a non-neighbour in $D_i$.}

 Suppose not. Note that since $C_i\neq \emptyset$, every vertex in $Z_i$ has a non-neighbour in $C_i$. Let $D_i\subseteq C_i$ be minimal such that every vertex in $Z_i$ has a non-neighbour in $D_i$. Then $D_i\neq \emptyset$, and so we have $|D_i|\geq 2k^2$. From the minimality of $D_i$, it follows that for every vertex $u\in D_i$, there is a vertex $z_u\in Z_i$ which is non-adjacent to $u$ and adjacent to all other vertices in $D_i$. Also, since $|D_i|\geq 2k^2$ and $C_1,\ldots, C_k$ partition $V(G)$ into $k$ stable sets, it follows that there $2k$ pairwise distinct vertices $u_1,\ldots, u_{2k}\in D_i$ such that $z_{u_1},\ldots, z_{u_{2k}}\in C_{i'}$ for some $i'\in \{1,\ldots, k\}\setminus \{i\}$; in particular, $\{z_{u_1},\ldots, z_{u_{2k}}\}$ is a stable set in $G$. Now, for every $j\in \{1,\ldots, 2k\}$, assuming $B_j=C_i\setminus N(z_j)$, it follows from the definition of $Z_i$ that $|B_j|\geq |C_i|/k$. Also, for all distinct $j,j'\in \{1,\ldots, 2k\}$, we have $|B_j\cap B_{j'}|<t$, as otherwise we may choose $A\subseteq B_j\cap B_{j'}$ with $|A|=t$, and then $G[A\cup \{u_j,u_{j'},z_j,z_{j'}\}]$ is isomorphic to $L_t$, which is impossible. But now by Lemma~\ref{incex}, we have
 $$\displaystyle |C_i|\geq \left|\bigcup_{j=1}^{2k}B_j\right|\geq \sum_{j=1}^{2k}\left|B_j\right|-\sum_{1\leq j<j'\leq 2k}|B_j\cap B_{j'}|> 2k\left(\frac{|C_i|}{k}\right)-4k^2t,$$
 which in turn implies that $|C_i|< 4k^2t$, a contradiction. This proves \eqref{st:smallD}.

 \medskip

Next, we define $\mathcal{S}_i$, for each $i\in \{1,\ldots, k\}$, to be the set of all maximum stable sets $S$ of $G$ with $k|S\cap C_i|\geq |C_i|$. The following is immediate from this definition and the fact that $k\alpha(G)\geq \sum_{i=1}^k|C_i|=|G|$.

\sta{\label{st:1/k} Every maximum stable set of $G$ belongs to $\mathcal{S}_i$ for some $i\in \{1,\ldots, k\}$.}

Moreover, we have:

 \sta{\label{st:Dhits} Let $i\in \{1,\ldots, k\}$ such that $|C_i|\geq 4k^2t$. Let $D_i\subseteq C_i$ with $|D_i|<2k^2$ be as promised by \eqref{st:smallD}. Then for every $S\in \mathcal{S}_i$, we have $S\cap D_i\neq\emptyset$.}

 Suppose not. Let $S\in \mathcal{S}_i$ with $S\cap D_i=\emptyset$. Since $D_i\neq \emptyset$, it follows that $S\not \subset C_i$ (as otherwise $C_i$ is a stable set of $G$ of cardinality more than $\alpha(G)$). Also, for every $z\in S\setminus C_i$, we have $C_i\cap S\subseteq C_i\setminus N(z)$, which along with the definition of $\mathcal{S}_i$, implies that $k|C_i\setminus N(z)|\geq k|C_i\cap S|\geq |C_i|$. We conclude that $S\setminus C_i$ is a non-empty subset of $Z_i$. Now, let $z\in S\setminus C_i\subseteq Z_i$ be a vertex with as many neighbours in $D_i$ as possible; then by \eqref{st:smallD}, $z$ has a non-neighbour $u\in D_i$, as well. Also, since $S$ is a maximum stable set of $G$ disjoint from $D_i$, it follows that $u$ has a neighbour $z'\in S$. From this and the fact that $C_i$ is a stable set of $G$, it follows that $z'\in (S\setminus C_i)\setminus \{z\}$. By the choice of $z$, there must be a vertex $u'\in D_i$ which is adjacent to $z$ and non-adjacent to $z'$; it follows that $u$ and $u'$ are distinct. Observe that $G[\{u,u',z,z'\}]$ is isomorphic to $M_2$. Moreover, since $S\in \mathcal{S}_i$, we have $|S\cap C_i|\geq |C_i|/k\geq 4kt>t$, and so we may choose $A\subseteq S\cap C_i$ with $|A|=t$. Note that since both $C_i$ and $S$ are stable sets of $G$ and $S\cap D_i=\emptyset$, it follows that $A$ and $\{u,u',z,z'\}$ are disjoint and anticomplete. But then $G[A\cup \{u,u',z,z'\}]$ is isomorphic to $L_t$, a contradiction. This proves \eqref{st:Dhits}.

 \medskip

Let us now put everything together. For every $i\in \{1,\ldots, k\}$, define $W_i$ as follows. If $|C_i|<4k^2t$, then let $W_i=C_i$, and if $|C_i|\geq 4k^2t$, let $W_i=D_i$, where $D_i\subseteq C_i$ with $|D_i|<2k^2$ is as in \eqref{st:smallD}. It follows that $|W_i|<4k^2t$. Let $W=\bigcup_{i=1}^kW_i$. Observe that $t\leq 2^{t-1}$ (which can be proved by an elementary induction on $t$). Also, since $c>1$, we have $2^{t+1}\leq c^{t+1}$. This yields:
$$|W|\leq 4k^3t\leq 4tc^{6t+6}\leq c^{7t+7}.$$

It remains to show that $W$ is a hitting set for $G$. Let $S$ be a maximum stable set of $G$. By \eqref{st:1/k}, we have $S\in \mathcal{S}_i$ for some $i\in \{1,\ldots, k\}$; in particular, we have $|S\cap C_i|\geq |C_i|/k>0$. If $|C_i|<4k^2t$, then $S\cap W_i=S\cap C_i\neq \emptyset$. Otherwise, by \eqref{st:Dhits}, we have $S\cap W_i=S\cap D_i\neq \emptyset$. In conclusion, we have shown that $W_i$, and so $W$, intersects $S$. This completes the proof of Theorem~\ref{mainparallel}.
\end{proof}

Now Theorem~\ref{properP51} becomes immediate:

\begin{theorem}
    Let $H$ be a proper induced subgraph of $P_5$. Then every $H$-free graph $G$ satisfies $\eta(G)\leq \omega(G)^{14}$.
\end{theorem}
\begin{proof}
    It is easy to check that $H$ is isomorphic to an induced subgraph of $P_4$, $S_{(2,1)}$, or $L_1$. But then the result follows from Theorems~\ref{perfectthm}, \ref{mainKrs} and \ref{mainparallel}.
\end{proof}

\section{Acknowledgement}

Our thanks to the anonymous referees for suggesting a number of improvements, and to Kaiyang Lan for pointing out an error in the proof of \eqref{st:smallD} (which is now fixed in this version).

\end{document}